\documentclass{article}
\usepackage{snapshot}

\pdfoutput=1
\usepackage{amsthm}

\usepackage[affil-sl]{authblk}

\usepackage[%
  backend = biber,%
  style = numeric,%
  natbib = true,%,
  url = false,%
  doi = false,%
  eprint = false,%
  backref=true,%
]{biblatex}

\usepackage{kpfonts, sectsty}

\usepackage[onehalfspacing, nodisplayskipstretch]{setspace}

\usepackage{xcolor}

\usepackage{hyperref}

\newtheoremstyle{plain}
  {3mm}                         % Space above theorem and previous line.
  {3mm}                         % Space below theorem box and next line.
  {\slshape}                    % Use slanted font in body of theorem.
  {}                            % Indent amount from margin. (Here no
  {\color{Blue}\bfseries}       % Theorem head font.
  {.}                           % Punctation after theorem head.
  {.5em}                        % Space after theorem head.
  {}                            % Theorem head specification. 

\newtheoremstyle{definition}
  {3mm}                         % Space above theorem and previous line.
  {3mm}                         % Space below theorem box and next line.
  {}                            % Use slanted font in body of theorem.
  {}                            % Indent amount from margin. (Here no
  {\color{Blue}\bfseries}       % Theorem head font.
  {.}                           % Punctation after theorem head.
  {.5em}                        % Space after theorem head.
  {}                            % Theorem head specification

\theoremstyle{plain}
\newtheorem{Theorem}{Theorem}[section]
\newtheorem{Fact}[Theorem]{Fact}
\newtheorem{Lemma}[Theorem]{Lemma}
\newtheorem{Proposition}[Theorem]{Proposition}
\newtheorem{Corollary}[Theorem]{Corollary}

\newtheorem*{CentralSetsTheorem}{Central Sets Theorem}
\newtheorem*{HalesJewettTheorem}{Hales--Jewett Theorem}

\theoremstyle{definition}
\newtheorem{Definition}[Theorem]{Definition}
\newtheorem{Remark}[Theorem]{Remark}

\title{%
  \textcolor{Scarlet}{A new and simpler noncommutative central sets theorem}%
}

\author{%
  John H. Johnson Jr.\\
  \small
  \href{mailto:johnson.5316@osu.edu}{\url{johnson.5316@osu.edu}}
}
\affil{%
  \small%
  Department of Mathematics\\
  The Ohio State University\\
  Columbus, Ohio
}

\sectionfont{\color{Scarlet}}
\subsectionfont{\color{Blue}}
\definecolor{Scarlet}{rgb}{0.78, 0.11, 0.0}
\definecolor{Blue}{rgb}{0.0, 0.42, 0.47}
\definecolor{Green}{rgb}{0.39, 0.71 ,0.0}
\definecolor{Orange}{rgb}{0.608, 0.306, 0}

\newcommand{\allNonemptyFiniteSubsetsOf}[1]{\mathcal{P}_f \left( #1 \right)}

\newcommand{\Define}[1]{\textcolor{Orange}{\textsl{\textbf{#1}}}}

\newcommand{\bbN}{\mathbb{N}}

\newcommand{\calT}{\mathcal{T}}
\newcommand{\calI}{\mathcal{I}}

\newcommand{\CartesianProduct}{\varprod}
\newcommand{\Ceiling}[1]{\left\lceil #1 \right\rceil}
\newcommand{\Cset}{$C$~set}
\newcommand{\Csets}{$C$~sets}
\newcommand{\Closure}[1]{\overline{#1}}
\newcommand{\FiniteSubsets}{\mathcal{P}_f}
\newcommand{\Floor}[1]{\left\lfloor #1 \right\rfloor}
\newcommand{\NonUniformSplitting}{\mathcal{I}}
\newcommand{\UniformSplitting}{\mathcal{J}}
\newcommand{\Jset}{$J$~set}
\newcommand{\Jsets}{$J$~sets}

\newcommand{\suchThat}{|\;}

\newcommand{\la}{\langle}
\newcommand{\ra}{\rangle}

\newcommand{\aVariableLetter}{\star}

\newcommand{\TheoremRef}[1]{%
  \hyperref[theorem:#1]{%
    \textcolor{Blue}{Theorem~\ref{theorem:#1}}%
  }}

\newcommand{\CorollaryRef}[1]{%
  \hyperref[corollary:#1]{%
    \textcolor{Blue}{Corollary~\ref{corollary:#1}}%
  }}

\newcommand{\LemmaRef}[1]{%
  \hyperref[lemma:#1]{%
    \textcolor{Blue}{Lemma~\ref{lemma:#1}}%
  }}

\newcommand{\PropositionRef}[1]{%
  \hyperref[proposition:#1]{%
    \textcolor{Blue}{Proposition~\ref{proposition:#1}}%
  }}

\newcommand{\RemarkRef}[1]{%
  \hyperref[remark:#1]{%
    \textcolor{Blue}{Remark~\ref{remark:#1}}%
  }}

\newcommand{\DefinitionRef}[1]{%
  \hyperref[definition:#1]{%
    \textcolor{Blue}{Definition~\ref{definition:#1}}%
  }}

\newcommand{\SectionRef}[1]{%
  \hyperref[section:#1]{%
    \textcolor{Blue}{Section~\ref{section:#1}}%
  }}

\newcommand{\OutsideTheoremRef}[1]{%
  \textcolor{Blue}{#1}%
  }

\newcommand{\OutsideDefinitionRef}[1]{%
  \textcolor{Blue}{#1}%
  }

\newcommand{\OutsideSectionRef}[1]{%
  \textcolor{Blue}{#1}%
  }

\newcommand{\SetFunction}[2]{\hbox{${}^{\hbox{$#1$}}{\hskip -2 pt #2}$}}

\hypersetup{%
  citecolor = Orange,
  colorlinks,
  linkcolor = Orange,
  urlcolor = Scarlet%
}

\begin{document}
\maketitle

\begin{abstract}
  Using dynamics, Furstenberg defined the concept of a central subset of positive integers and proved several powerful combinatorial properties of central sets.  
  Later using the algebraic structure of the Stone--\v{C}ech compactification, Bergelson and Hindman, with the assistance of B. Weiss, generalized the notion of a central set to any semigroup and extended the most important combinatorial property of central sets to the central sets theorem.
  Currently the most powerful formulation of the central sets theorem is due to De, Hindman, and Strauss in \cite[\OutsideTheoremRef{Corollary 3.10}]{De:2008:NewStrongerCST}.
  However their formulation of the central sets theorem for noncommutative semigroups is, compared to their formulation for commutative semigroups, complicated.
  In this paper I prove a simpler (but still equally strong) version of the noncommutative central sets theorem in \CorollaryRef{central-sets-theorem}.
\end{abstract}

\section{Introduction}
\label{section:introduction}
Furstenberg, in his book connecting dynamics to Ramsey theory, defined the concept of a \emph{central subset of positive integers} \cite[\OutsideDefinitionRef{Definition 8.3}]{Furstenberg:1981:RecurrenceNumThy} and proved several important properties of central sets using notions from topological dynamics.
One such property is that whenever a central set is finitely partitioned, at least one cell of the partition contains a central set \cite[\OutsideTheoremRef{Theorem 8.8}]{Furstenberg:1981:RecurrenceNumThy}.
Most of the remaining important properties of central sets are derivable from a powerful combinatorial theorem \cite[\OutsideTheoremRef{Proposition 8.21}]{Furstenberg:1981:RecurrenceNumThy}.
(A bit later, in \TheoremRef{commutative-central-sets-theorem} on page \pageref{theorem:commutative-central-sets-theorem}, I state the most powerful formulation of \cite[\OutsideTheoremRef{Proposition 8.21}]{Furstenberg:1981:RecurrenceNumThy}.)
Furstenberg used his combinatorial theorem to prove Rado's theorem (see the sufficiency condition of \cite[\OutsideTheoremRef{Theorem 5} on page 74]{Graham:1990:RamseyTheory} or Rado's published dissertation \cite[\OutsideTheoremRef{Satz IV} on page 445]{Rado:1933:StudiesOnComb}) by showing that given any central set and any $m \times n$ matrix $M$, with integer entries, that satisfies the ``columns condition'' we can find a vector $\vec{x}$ all of whose entries are in the central set with $M\vec{x} = 0$.

Inspired by the fruitful interaction between ultrafilters on semigroups and Ramsey theory, Bergelson and Hindman, with the assistance of B.~Weiss, later proved an algebraic characterization of a central subset of positive integers \cite[\OutsideSectionRef{Section 6}]{Bergelson:1990:DynRamseyTheory}.
This algebraic definition has several advantages over the original dynamical definition: one advantage is that the algebraic definition is simple and easily generalizes to any semigroup. 
\begin{Definition}
  \label{definition:central-set}
   Let $(S, \cdot)$ be a discrete semigroup and let $A \subseteq S$.
   Then $A$ is a \Define{central set} if and only if there exists an idempotent $p$ in the smallest ideal of $\beta S$ with $A \in p$.
\end{Definition}

Another advantage of the algebraic definition is that central sets  are ``partition regular''~---~that is, in any finite partition of a central set at least one cell of the partition is a central set \cite[\OutsideTheoremRef{Theorem 8.8}]{Furstenberg:1981:RecurrenceNumThy}~---~is, because of standard properties of ultrafilters, immediate from the definition.
More importantly, Furstenberg's combinatorial theorem \cite[\OutsideTheoremRef{Proposition 8.21}]{Furstenberg:1981:RecurrenceNumThy} follows from a relatively simple recursive construction.

\begin{Remark}
  \label{remark:dynamical-definition-central-set}
  The dynamical definition of a central set also extends naturally to an arbitrary semigroup; the fact that the algebraic and dynamical definitions are equivalent was proved by H.~Shi and H.~Yang in \cite{Shi:1996:DynCharCentralSets}.
  Besides this introduction I will usually not emphasize the dynamical point-of-view.
\end{Remark}

Using the algebraic structure of $\beta S$, the Stone--\v{C}ech compactification of a discrete semigroup $S$, and a more sophisticated recursive construction De, Hindman, and Strauss proved the (currently) strongest version of the central sets theorem in \cite{De:2008:NewStrongerCST}.
(The central sets theorem is what we shall call the main combinatorial property central sets satisfy.)
We first state the central sets theorem for commutative semigroups: the statement of the noncommutative version is more complicated and forms the main focus of this paper.

In the statement of the commutative central sets theorem, and in the remainder of this paper, we let $\allNonemptyFiniteSubsetsOf{X}$ denote the collection of all nonempty finite subsets of a (typically nonempty) set $X$, and we let $\SetFunction{A}{B}$ denote the collection of all functions with domain $A$ and codomain $B$.

\begin{Theorem}[Commutative central sets theorem]
  \label{theorem:commutative-central-sets-theorem}
  Let $(S, +)$ be a commutative semigroup and let $A \subseteq S$ be central.
  For typographical convenience we let $\calT = \SetFunction{\bbN}{S}$.
  Then there exist functions $\alpha \colon \allNonemptyFiniteSubsetsOf{\calT} \to S$ and $H \colon \allNonemptyFiniteSubsetsOf{\calT} \to \allNonemptyFiniteSubsetsOf{\bbN}$ that satisfy the following two statements:
  \begin{itemize}
    \item[(1)] If $F$ and $G$ are both in $\allNonemptyFiniteSubsetsOf{\calT}$ with $F \subsetneq G$, then $\max H(F) < \min H(G)$.

    \item[(2)] Whenever $m$ is a positive integer, $G_1, G_2, \ldots, G_m$ is a finite sequence in $\allNonemptyFiniteSubsetsOf{\calT}$ with $G_1 \subsetneq G_2 \subsetneq \cdots \subsetneq G_m$ and for every $i \in \{1, 2, \ldots, m\}$ we have $f_i \in G_i$, then
        $
          \sum_{i=1}^m \bigl( \alpha(G_i) + \sum_{t \in H(G_i)} f_i(t) \bigr) \in A.
        $
    \end{itemize}
\end{Theorem}
\begin{proof}
  This was proved by De, Hindman, and Strauss as \cite[\OutsideTheoremRef{Theorem 2.2}]{De:2008:NewStrongerCST}.
\end{proof}

In the same paper, De, Hindman, and Strauss also formulated and proved a strong version of the central sets theorem for arbitrary semigroups \cite[\OutsideTheoremRef{Corollary 3.10}]{De:2008:NewStrongerCST}.
The statement of this version of the central sets theorem is necessarily complicated because the underlying semigroup may be noncommutative.
(When the underlying semigroup is commutative, \cite[\OutsideTheoremRef{Corollary 3.10}]{De:2008:NewStrongerCST} reduces to \TheoremRef{commutative-central-sets-theorem}.)
In the case of noncommutativity it is usually \emph{not sufficient}~---~for both combinatorial (see \cite[\OutsideSectionRef{Section 1}]{Bergelson:1992:RamseyThySemigrps}) and algebraic (see \cite[\OutsideTheoremRef{Theorem 1.13}]{Johnson:2011:Dissertation} reasons~---~to simply perform the obvious translation of \TheoremRef{commutative-central-sets-theorem} to an arbitrary semigroup. Roughly speaking, the proper translation requires splitting up each translate $\alpha(G_i)$ into several parts.

To better explain this difference in formulation let's first consider a special case of \TheoremRef{commutative-central-sets-theorem}.
\begin{Corollary}
  \label{corollary:central-implies-Jset}
  Let $(S, +)$ be a commutative semigroup and $A \subseteq S$ central.
  Then for every $F \in \SetFunction{\bbN}{S}$ there exist $a \in S$ and $H \in \allNonemptyFiniteSubsetsOf{\bbN}$ such that for every $f \in F$ we have
    $
      a + \sum_{t \in H} f(t) \in A.
    $
\end{Corollary}
\begin{proof}
  Pick functions $\alpha$ and $H$ as guaranteed by \TheoremRef{commutative-central-sets-theorem}.
  Let $F \in \SetFunction{\bbN}{S}$, put $m = 1$, and observe from conclusion (2) of \TheoremRef{commutative-central-sets-theorem} we have that for every $f \in F$
    \[
      \alpha(F) + \sum_{t \in H(F)} f(t) \in A.
    \]
\end{proof}

To perform the appropriate splitting of the ``$a$'' in \CorollaryRef{central-implies-Jset} (or ``$\alpha(G_i)$'' in \TheoremRef{commutative-central-sets-theorem}) De, Hindman, and Strauss introduced, via \cite[\OutsideDefinitionRef{Definition~3.1}]{De:2008:NewStrongerCST}, the notation $\NonUniformSplitting_m$:

\begin{Definition}
  \label{definition:nonuniform-splitting}
  For every positive integer $m$ define
  \begin{align*}
    \NonUniformSplitting_m &= \bigl\{  \bigl(H(1), H(2), \ldots, H(m)\bigr) \bigm\suchThat \mbox{every $H(i) \in \allNonemptyFiniteSubsetsOf{\bbN}$ and} \\
      &\hspace{2em}\mbox{for each $i \in \{1, 2, \ldots, m-1\}$ we have $\max H(i) < \min H(i+1)$} \bigr\}. 
    \end{align*}
\end{Definition}

Using this notation and their version of the Central Sets Theorem \cite[\OutsideTheoremRef{Corollary 3.10}]{De:2008:NewStrongerCST} it's possible to, by analogy with \CorollaryRef{central-implies-Jset}, derive as a special case the following fact:
\begin{Proposition}
  \label{proposition:central-implies-old-Jset}
  Let $(S, \cdot)$ be a semigroup and let $A \subseteq S$ be central.
      Then for every $F \in \SetFunction{\bbN}{S}$ there exists a positive integer $m$ such that for some $a \in S^{m+1}$ and $H \in \NonUniformSplitting_m$ and for every $f \in F$ we have
    $
      \prod_{i=1}^m\bigl( a(i) \cdot \prod_{t \in H(i)} f(t) \bigr)a(m+1) \in A.
    $
\end{Proposition}
\begin{Remark}
  \label{remark:small-example-on-nonuniform-notation}
  To understand this last expression it may be helpful to consider a small representative example.
  For instance, suppose $m = 3$, $H(1) = \{2, 5\}$, $H(2) = \{10\}$, and $H(3) = \{21, 22, 50\}$, then
  \begin{align*}
    \Bigl( \prod_{i=1}^3 \bigl( a(i) \cdot \prod_{t \in H(i)} f(t) \bigr) \Bigr)\cdot a(4) &= a(1) \cdot f(2)f(5) \cdot a(2) \\
    &\hspace{2em}\cdot f(10) \cdot a(3) \cdot f(21)f(22)f(50) \cdot a(4).
  \end{align*}
\end{Remark}

Note that when the underlying semigroup is commutative, by appropriate rewriting, \PropositionRef{central-implies-old-Jset} reduces to \CorollaryRef{central-implies-Jset}.

In De's, Hindman's, and Strauss's formulation and their proof~---~and as illustrated in the small example in \RemarkRef{small-example-on-nonuniform-notation}~---~the choice of $H \in \NonUniformSplitting_m$ may be ``nonuniform'', that is, each $H(i)$ may have a different size.
My main result in this paper is to show that, without loss of generality, we  \emph{can} require that each $H(i)$ be a singleton.
From this observation we can produce a new and simpler formulation of the noncommutative central sets theorem.
As an easily stated special case of this theorem we have
\begin{Proposition}
  \label{proposition:central-implies-new-Jset}
  Let $(S, \cdot)$ be a semigroup and let $A \subseteq S$ be central.
  Then for every $F \in \SetFunction{\bbN}{S}$ there exist $m \in \bbN$, $a \in S^{m+1}$, and a strictly increasing sequence $\la t_i \ra_{i=1}^m$ of positive integers such that for every $f \in F$ we have
    $
    \Bigl( \prod_{i=1}^m\bigl(a(i) f(t_i)\bigr) \Bigr) a(m+1) \in A.
    $
\end{Proposition}

In \SectionRef{properties-of-Jsets}, we will name those sets that satisfy the conlusion of \PropositionRef{central-implies-new-Jset} and prove some important combinatorial and algebraic properties of these sets; in \SectionRef{simpler-central-sets-theorem}, we state and prove our new and simpler formulation of the noncommutative central sets theorem; and, in the final \SectionRef{equivalence-of-central-set-theorems}, we prove that our formulation of the central sets theorem, while simpler, is equivalent to \cite[\OutsideTheoremRef{Corollary 3.10}]{De:2008:NewStrongerCST}.

I'll end this introduction by presenting a brief review of the algebraic structure of the Stone--\v{C}ech compactification of a discrete semigroup, with references to proofs of these facts.

Given a discrete nonempty space $S$ we take the points of $\beta S$ to be the collection of all ultrafilters on $S$.
We identify points of $S$ with the principal ultrafilters in $\beta S$.
(Thus we pretend that $S \subseteq \beta S$.)
Given $A \subseteq S$ put $\overline{A} = \{ p \in \beta S \mid A \in p\}$.
Then the collection $\{ \overline{A} \mid A \subseteq S \}$ is a basis for a compact Hausdorff topology on $\beta S$.
This topology is the Stone--\v{C}ech compactification of the discrete space $S$.
The proofs for all of these assertions can be found in \cite[\OutsideSectionRef{Sections 3.2 and 3.3}]{Hindman:1998}.

Given a discrete semigroup $(S, \cdot)$, we can extend the semigroup operation to $\beta S$ \cite[\OutsideTheoremRef{Theorem 4.1}]{Hindman:1998} such that for $p$ and $q$ in $\beta S$ and every $A \subseteq S$ we have $A \in p \cdot q$ if and only if $\{ x \in S \mid x^{-1}A \in q \} \in p$ \cite[\OutsideTheoremRef{Theorem 4.12}]{Hindman:1998} where $x^{-1}A = \{ y \in S \mid xy \in A\}$.
With this operation $(\beta S, \cdot)$ becomes a compact Hausdorff right-topological semigroup: \Define{right-topological} means that for every $q \in \beta S$, the function $p \mapsto p \cdot q$, whose domain and codomain are both $\beta S$, is continuous.

\begin{Fact}
  \label{fact:crts-has-idempotent-and-smallest-ideal}
  Let $T$ be a compact Hausdorff right-topological semigroup.
  \begin{itemize}
    \item[(a)] $T$ contains at least one idempotent, that is, there exists $x \in T$ such that $x = x \cdot x$.

    \item[(b)] $T$ contains an ideal, called the smallest ideal and denoted by $K(T)$, that is contained in every ideal of $T$.  Additionally, $K(T)$ also contains at least one idempotent.
    \end{itemize}
\end{Fact}
\begin{proof}
  The proof of statements \textsl{(a)} and \textsl{(b)} are given in \cite[\OutsideTheoremRef{Theorems 2.5 and 2.8}]{Hindman:1998}, respectively.
  (The proof of statement \textsl{(a)} is originally due to Ellis in \cite[\OutsideTheoremRef{Lemma 1}]{Ellis:1958:DistalTransformationGroups}.)
\end{proof}

\subsection*{Acknowledgements}
\label{section:ackonwledgements}
This new and simpler noncommutative central sets theorem is part of my dissertation research conducted under the guidance of Neil Hindman; I thank Dr.~Hindman for his excellent advisement.
  Also, I thank both Vitaly Bergelson and Neil Hindman for recently encouraging me to write up this result for publication, their extraordinary patience as this paper was (slowly) being written and revised, and reading and commenting on an early draft of this paper.

\section{Simpler definition of {\Jsets} and its algebraic and combinatorial properties}
\label{section:properties-of-Jsets}
Our first task will be to define and study the properties of {\Jsets}~---~these are the sets that satisfy the conclusion of \PropositionRef{central-implies-new-Jset}.
The definition we give is new and simpler than the one given in \cite[\OutsideDefinitionRef{Definition 2.3(d)}]{Hindman:2010:CartesianProductCsets}; and, this simplicity is the core of the simpler formulation of the central sets theorem.
Finally, note that I will define a bit more notation (such as `$\UniformSplitting_m$', `$\calT$', and `$x(m, a, t, f)$') then strictly necessary to define a {\Jset}: this extra notation is useful when formulating the central sets theorem later.

\begin{Definition}
  \label{definition:new-Jset}
  Let $(S, \cdot)$ be a semigroup.
  \begin{itemize}
    \item[(a)]
      For each positive integer $m$ put
      $
        \UniformSplitting_m = \{ (t_1, t_2, \ldots, t_m) \in \bbN^{m} \mid t_1 < t_2 < \cdots < t_m \}.
      $

   \item[(b)]
     Put $\calT(S) = \SetFunction{\bbN}{S}$, and if the semigroup is clear from context we will instead write $\calT$ for $\calT(S)$.

   \item[(c)]
     For all positive integers $m$, every $a \in S^{m+1}$, every $t \in \UniformSplitting_m$, and for all $f \in \calT$, put
      \[
        x(m, a, t, f) = \Bigl( \prod_{i=1}^m \bigl( a(i) f(t_i) \bigr) \Bigr) a(m+1).
      \]

    \item[(d)]
      We call $A \subseteq S$ a \Define{{\Jset} (in $S$)} if and only if for every $F \in \allNonemptyFiniteSubsetsOf{\calT}$ there exist $m \in \bbN$, $a \in S^{m+1}$, and $t \in \UniformSplitting_m$ such that for every $f \in F$ we have $x(m, a, t, f) \in A$.
  \end{itemize}
\end{Definition}
\begin{Remark}
  \label{remark:usual-warning-on-Jsets}
  I must point out that {\Jsets} are \emph{not} named after the author!
  The term {\Jsets} is derived from an earlier term of ``$J_Y$~set'' introduced in \cite[\OutsideDefinitionRef{Definition 2.4(b)}]{Hindman:1996:CentralSetsCombChar}.
\end{Remark}

We start with an easy observation showing that we have some control in dictating the position of the starting $t_1$ term in $x(m, a, t, f)$.
(We won't use this observation in this section but it will immediately be used in \SectionRef{simpler-central-sets-theorem} for \TheoremRef{idempotents-and-Csets}.)

\begin{Lemma}
  \label{lemma:starting-position}
  Let $(S, \cdot)$ be a semigroup, $A \subseteq S$ a {\Jset}, and let $N$ be a positive integer.
  For every $F \in \allNonemptyFiniteSubsetsOf{\calT}$ there exist $m \in \bbN$, $a \in S^{m + 1}$, and $t \in \UniformSplitting_m$ with $t_1 > N$ such that for every $f \in F$ we have $x(m, a, t, f) \in A$.
\end{Lemma}
\begin{proof}
  For every $f \in F$ define $g_f \in \calT$ by $g_f(t) = f(t + N)$.
  Pick $m \in \bbN$, $a \in S^{m+1}$, and $u \in \UniformSplitting_m$ such that for every $f \in F$ we have $x(m, a, u, g_f) \in A$.
  Define $t \in \UniformSplitting_m$ by $t_i = u_i + N$ for every $i \in \{1, 2, \ldots, m \}$, then $x(m, a, t, f) = x(m, a, u, g_f) \in A$.
\end{proof}

Our next task, which is more substantial, is to show that {\Jsets} are partition regular, that is, whenever a {\Jset} is finitely partitioned at least one cell in the partition will be a {\Jset}.
We start with this result since, although our proof is based on \cite[\OutsideTheoremRef{Theorem 2.14}]{Hindman:2010:CartesianProductCsets}, it nicely illustrates how proofs can be simplified using the simpler definition of a {\Jset}.
(For instance compare the computation of the translate $a \in S^{m + 1}$ (their $c \in S^{m + 1}$) in \cite[\OutsideTheoremRef{Theorem 2.14}]{Hindman:2010:CartesianProductCsets} with my computation of $a$.)
Like \cite[\OutsideTheoremRef{Theorem 2.14}]{Hindman:2010:CartesianProductCsets} I will use the powerful and flexible Hales--Jewett theorem \cite[\OutsideTheoremRef{Theorem 1}]{Hales:1963:Regularity}.
Accordingly, I'll state some needed definitions and the version of the Hales--Jewett theorem I will use.
\begin{Definition}
  \label{definition:variable-word-and-combinatorial-line}
  Let $k$ and $N$ both be positive integers and let $\aVariableLetter$ be a symbol not in $\{ 1, 2, \ldots, k \}$.
  \begin{itemize}
    \item[(a)]
      We call $w(\aVariableLetter) \in (\{ 1, 2, \ldots, k \} \cup \{ \aVariableLetter \})^N$ a \Define{variable word of length $N$ over (the alphabet) $\{ 1, 2, \ldots, k \}$} if and only if $\aVariableLetter$ occurs at least once in the coordinates of $w(\aVariableLetter)$.

    \item[(b)]
      Let $w(\aVariableLetter)$ be a variable word of length $N$ over $\{ 1, 2, \ldots, k \}$.
      We call \[\bigl\{ w(i) \suchThat i \in \{ 1, 2, \ldots k \} \bigr\}\] a \Define{combinatorial line}, where $w(i)$ is the element in  $\{ 1, 2, \ldots, k \}^N$ formed when we replace each occurrence of $\aVariableLetter$ in $w(\aVariableLetter)$ by $i$
  \end{itemize}
\end{Definition}

\begin{HalesJewettTheorem}
  \label{theorem:hales-jewett-theorem}
  For every positive integer $k$ there exists a positive integer $N$ such that if $c \colon \{ 1, 2, \ldots, k \}^N \to \{ 1, 2 \}$, then there exists a combinatorial line on which $c$ is constant.
\end{HalesJewettTheorem}

Observe that to prove that {\Jsets} are partition regular it suffices to prove the following theorem.

\begin{Theorem}
  \label{theorem:Jsets-are-partition-regular}
  Let $(S, \cdot)$ be a semigroup and let $A$ and $B$ both be subsets of $S$.
  If $A \cup B$ is a {\Jset}, then either $A$ is a {\Jset} or $B$ is a {\Jset}.
\end{Theorem}
\begin{proof}
  If $A$ is a {\Jset}, then we are obviously done.
  So, we may assume that $A$ is not a {\Jset}.
  Under this additional assumption we show, via the Hales--Jewett theorem and (somewhat tedious) rewriting, that $B$ is a {\Jset}.

  Let $F \in \allNonemptyFiniteSubsetsOf{\calT}$.
  Since $A$ is not a {\Jset} we may pick $G \in \allNonemptyFiniteSubsetsOf{\calT}$ such that for every $m \in \bbN$, every $a \in S^{m+1}$, and every $t \in \UniformSplitting_m$ there exists $g \in G$ such that $x(m, a, t, g) \not\in A$.
  Put $H = F \cup G$, put $k = |H|$, and enumerate $H$ as $\{ h_1, h_2, \ldots, h_k \}$.
  Finally pick a positive integer $N$ as guaranteed by the Hales--Jewett theorem for $k$.

  Fix $d$ in $S$ and for every $w = (w_1, w_2, \ldots, w_N) \in \{ 1, 2, \ldots, k \}^N$ define $g_w \in \calT$ by
  \[
    g_w(t) = \prod_{i = 1}^N \bigl(d \cdot h_{w_i} (N t + i) \bigr).
  \]
  Since, by hypothesis, $A \cup B$ is a {\Jset} and $\bigl\{ g_w \bigm| w \in \{ 1, 2, \ldots, k \}^N \bigr\} \in \allNonemptyFiniteSubsetsOf{\calT}$ we may pick $p \in \bbN$, $b \in S^{p+1}$, and $s \in \UniformSplitting_p$ such that for every $w \in \{ 1, 2, \ldots, k \}^N$ we have $x(p, b, s, g_w) \in A \cup B$.

  Now define the function $c \colon \{ 1, 2, \ldots, k \}^N \to \{ 1, 2 \}$ by
  \[
    c(w) = 
    \begin{cases}
      1 & \mbox{if $x(p, b, s, g_w) \in A$, and} \\
      2 & \mbox{otherwise.}
    \end{cases}
  \]
  By the Hales--Jewett theorem we may pick a variable word $w(\aVariableLetter)$ over $\{ 1, 2, \ldots, k \}$ such that the combinatorial line $\bigl\{ w(i) \bigm\suchThat i \in \{ 1, 2, \ldots, k\} \bigr\}$ is constant on $c$.
  
  We now claim that there exist $m \in \bbN$, $a \in S^{m+1}$, and $t \in \UniformSplitting_m$ such that for every $i \in \{ 1, 2, \ldots, k \}$ we have $x(m, a, t, h_i) = x(p, b, s, g_{w(i)})$.
  If this claim is true, then by our assumption that $A$ is not a {\Jset} and by our choice of $N$ via the Hales--Jewett theorem it follows that, in particular, $x(m, a, t, f) \in B$ for all $f \in F$.
  Hence $B$ is a {\Jset}, under the additional assumption of our claim.

  Therefore we now spend the rest of our proof showing this claim.
  To this end we need a bit more notation and an additional assumption on our semigroup.
  Let $r$ be the number of occurrences of the variable letter in $w(\aVariableLetter)$.
  For each $i \in \{ 1, 2, \ldots, r \}$ let $v_i$ denote the position of the $i$th occurrence of the variable letter in $w(\aVariableLetter)$.
  Define
  \begin{align*}
    \Lambda_1 &= \bigl\{ i \in \{ 1, 2, \ldots, N \} \bigm\suchThat i < v_1 \bigr\}, \\
    \Lambda_{r + 1} &= \bigl\{ i \in \{ 1, 2, \ldots, N \} \bigm\suchThat v_r < i \bigr\}, \mbox{ and} \\
  \end{align*}
  for each $j \in \{2, 3, \ldots, r\}$ define
  \begin{align*}
    \Lambda_j &= \bigl\{ i \in \{ 1, 2, \ldots, N \} \bigm\suchThat v_{j - 1} < i < v_j \bigr\}.
  \end{align*}
  (Hence the $\Lambda_j$'s group the adjacent non-variable letters in $w(\aVariableLetter)$.)
  Since we shall be performing modulo arithmetic (with modulus $r$) on the index of the $v_i$'s and $\Lambda_j$'s we also put $v_0 = v_r$ and $\Lambda_0 = \Lambda_{r}$.
  We also assume that our semigroup $S$ has a two-sided identity $1$; if our semigroup doesn't have an identity, we can just append an identity to our old semigroup.
  We do this because in what follows we may consider empty products, and we shall interpret an empty product as multiplication by $1$.
  With this notational setup done we can now give $m \in \bbN$, $a \in S^{m+1}$, and $t \in \UniformSplitting_m$ that works to prove our claim.

  Put $m = p \cdot r$, define $a \in S^{m + 1}$ for every $i \in \{ 1, 2, \ldots, m + 1 \}$ by
  \[
    a(i) =
    \begin{cases}
      b(1) \cdot \Bigl( \prod_{j \in \Lambda_1} \bigl( d \cdot h_{w_j}(N \cdot s_1 + j) \bigr) \Bigr) \cdot d & \mbox{if $i = 1$,} \\
      & \\
      \Bigl( \prod_{j \in \Lambda_{i \bmod{r}}} \bigl( d \cdot h_{w_j}(N \cdot s_{\Ceiling{i/r}} + j) \bigr) \Bigr) \cdot d & \mbox{if $i \not\equiv 1 \pmod{r}$,} \\
      & \\
      \Bigl( \prod_{j \in \Lambda_{r + 1}} \bigl( d \cdot h_{w_j}(N \cdot s_{\Floor{i/r}} + j) \bigr) \Bigr) \cdot b(\Ceiling{i/r}) \\
      \hspace{1em} \cdot \Bigl( \prod_{i \in \Lambda_1} \bigl( d \cdot h_{w_j}(N \cdot s_{\Ceiling{i/r}} + j) \bigr) \Bigr) \cdot d &\mbox{if $i \ne 1$, $i \ne m + 1$,} \\
      & \hspace{2em}\mbox{and $i \equiv 1 \pmod{r}$, and} \\
      \Bigl( \prod_{j \in \Lambda_{r+1}} \bigl( d \cdot h_{w_j}(N \cdot s_{\Floor{i/r}} + j) \bigr) \Bigr) \cdot b(p + 1) & \mbox{if $i = m + 1$,}
    \end{cases}
  \]
  and we define $t \in \UniformSplitting_m$ for every $i \in \{ 1, 2, \ldots, m \}$ by $t_i = N \cdot s_{\Ceiling{i/r}} + v_{i \bmod{r}}$, observing that $s_1 < s_2 < \cdots < s_p$ implies $t_1 < t_2 < \cdots < t_m$.
  Then for every $i \in \{ 1, 2, \ldots, k \}$ we have $x(m, a, t, h_i) = x(p, b, s, g_{w(i)})$.
  This proves our claim and so completes the proof of this theorem.

    (To see how this complicated formula was derived it may be helpful to consider a representative example.
    Suppose $k = 2$, $N = 6$, and $w(\aVariableLetter) = (\aVariableLetter, 1, 2, \aVariableLetter, \aVariableLetter, 1)$.
    Then $r = 3$; $v_1 = 1$, $v_2 = 4$, and $v_3 = 5$ (which is also $v_0$); and we have $\Lambda_1 = \emptyset$, $\Lambda_2 = \{2, 3\}$, $\Lambda_3 = \emptyset$ (which is also $\Lambda_0$), and $\Lambda_4 = \{ 6 \}$.
    Further, suppose $p = 2$, then $x(2, b, s, g_{w(\aVariableLetter)}) = b(1) g_{w(\aVariableLetter)}(s_1) b(2) g_{w(\aVariableLetter)}(s_2) b(3)$.
    Then $m = 2 \cdot 3$ and by the above formulas
    \begin{align*}
      a(1) &= b(1) \cdot d; a(2) = d \cdot h_1(6s_1 + 2) \cdot d \cdot h_2(6s_1 + 3) \cdot d; a(3) = d;\\
      a(4) &= d \cdot h_1(6s_1+6) \cdot b(2) \cdot d; a(5) =  d \cdot h_1(6s_2 + 2) \cdot d \cdot h_2(6s_2 + 3) \cdot d; a(6) = d;\\
      a(7) &= d \cdot h_1(6s_2 + 6) \cdot b(3);
    \end{align*}
    and $t_1 = 6\cdot s_1 + 1$; $t_2 = 6\cdot s_1 + 4$; $t_3 = 6\cdot s_1 + 5$; $t_4 =6\cdot s_2 + 1$; $t_5 = 6\cdot s_2 + 4$; and $t_6 = 6\cdot s_2 + 4$.
    A bit of checking shows that $x(2, b, s, g_{w(i)}) = x(7, a, t, h_i)$.)
\end{proof}

\TheoremRef{Jsets-are-partition-regular} allows us to quickly derive some algebraic information about {\Jsets}.

\begin{Definition}
  \label{definition:ideal-of-J(S)}
  If $(S, \cdot)$ is a semigroup, then define
  \[
    J(S) = \{ p \in \beta S \mid \mbox{for every $A \in p$ we have that $A$ is a {\Jset}} \}.
  \]
\end{Definition}

\begin{Corollary}
  \label{corollary:Jsets-and-the-ideal-J(S)}
  Let $(S, \cdot)$ be a semigroup and let $A \subseteq S$.
  Then $A$ is a {\Jset} if and only if $\Closure{A} \cap J(S) \ne \emptyset$.
\end{Corollary}
\begin{proof}
  \textcolor{Blue}{($\Rightarrow$)} By \cite[\OutsideTheoremRef{Theorem 3.11}]{Hindman:1998} and \TheoremRef{Jsets-are-partition-regular} we can pick an ultrafilter $p$ on $S$ such that $A \in p$ and every element of $p$ is a {\Jset}.
  
  \textcolor{Blue}{($\Leftarrow$)} Trivially true by definition of $J(S)$.
\end{proof}

\begin{Theorem}
  \label{theorem:ideal-of-J(S)}
  Let $(S, \cdot)$ be a semigroup.
  $J(S)$ is a nonempty closed (two-sided) ideal of $\beta S$.
\end{Theorem}
\begin{proof}
    In particular, $S$ is a {\Jset} and so by \CorollaryRef{Jsets-and-the-ideal-J(S)} we have $\Closure{S} \cap J(S) \ne \emptyset$.

    To see that $J(S)$ is closed let $p \not\in J(S)$ and pick $A \in p$ such that $A$ is not a {\Jset}.
    By definition of $J(S)$ we have $\overline{A} \cap J(S) = \emptyset$.
    Since $\Closure{A}$ is a (basic) open neighborhood of $p$, it follows that $J(S)$ is topologicaly closed in $\beta S$.

    Now let $p \in J(S)$ and $q \in \beta S$.
    To see that $J(S)$ is an ideal, we show that $p \cdot q \in J(S)$ and $q \cdot p \in J(S)$.

    We first show that $p \cdot q \in J(S)$.
    Let $A \in p \cdot q$, let $F \in \allNonemptyFiniteSubsetsOf{\calT}$, and put $B = \{ x \in S \mid x^{-1}A \in q \}$.
    Now $A \in p \cdot q$ if and only if $B \in p$ and hence $B$ is a $J$ set.
    Pick $m \in \bbN$, $b \in S^{m+1}$, and $t \in \UniformSplitting_m$ as guaranteed for $B$ with respect to the collection $F$.
    By definition of $B$ for every $f \in F$ we have $x(m, b, t, f)^{-1}A \in q$ and so
    $\bigcap_{f \in F} x(m, b, t, f)^{-1}A \ne \emptyset$.
    Pick $c \in \bigcap_{f \in F} x(m, b, t, f)^{-1}A$ and define $a \in S^{m+1}$ by 
    \[
    a(i) =
    \begin{cases}
      b(i) & \mbox{if $i \in \{1, 2, \ldots, m\}$} \\
      b(m+1)c & \mbox{if $i = m+1$}.
    \end{cases}
    \]
    Therefore for each $f \in F$ we have $x(m, a, t, f) \in A$.
    Hence $A$ is a {\Jset} and we have $p \cdot q \in J(S)$.

    To show that $q \cdot p \in J(S)$ is similar.
    Let $A \in q \cdot p$, let $F \in \allNonemptyFiniteSubsetsOf{\calT}$, and put $B = \{ x \in S \mid x^{-1}A \in p \}$.
    Now $A \in q \cdot p$ if and only if $B \in q$.
    Since $B \ne \emptyset$, pick $b \in B$ and so $b^{-1}A \in p$.
    Since $b^{-1}A$ is a $J$ set pick $m \in \bbN$, $c \in S^{m+1}$, and $t \in \UniformSplitting_m$ such that for all $f \in F$ we have $x(m, c, t, f) \in b^{-1}A$.
    Define $a \in S^{m+1}$ by
    \[
      a(i) =
      \begin{cases}
        b\cdot c(1) & \mbox{if $i = 1$}\\
        c(i) & \mbox{if $i \in \{2, 3, \ldots, m+1\}$}.
      \end{cases}
    \]
    Therefore for each $f \in F$ we have $x(m, a, t, f) \in A$.
    Hence $A$ is a {\Jset} and we have $q \cdot p \in J(S)$.
\end{proof}
\begin{Corollary}
  \label{corollary:J(S)-contains-smallest-ideal}
  Let $(S, \cdot)$ be a semigroup.
  Then $K( \beta S ) \subseteq J(S)$.
\end{Corollary}
\begin{proof}
  By \TheoremRef{ideal-of-J(S)} $J(S)$ is an ideal and the smallest ideal $K( \beta S )$ is contained in every ideal of $\beta S$.
\end{proof}

I'll end this section by showing that we can provide a ``purely'' algebraic proof of \CorollaryRef{J(S)-contains-smallest-ideal} by exploiting the relationship between piecewise syndetic sets and the smallest ideal \cite[\OutsideTheoremRef{Theorems 4.39 and 4.40}]{Hindman:1998}.

\begin{Definition}
  \label{definition:piecewise-syndetic}
  Let $(S, \cdot)$ be a semigroup and $A \subseteq S$.
  We call $A$ a \Define{piecewise syndetic set} if and only if there exists $G \in \allNonemptyFiniteSubsetsOf{S}$ such that for all $F \in \allNonemptyFiniteSubsetsOf{S}$ there exists $x \in S$ with $Fx \subseteq \bigcup_{t \in G} t^{-1}A$.
\end{Definition}

\begin{Theorem}
  \label{theorem:piecewise-syndetic-implies-Jset}
  Let $(S, \cdot)$ be a semigroup and $A \subseteq S$.
  If $A$ is piecewise syndetic, then $A$ is a {\Jset}.
\end{Theorem}
\begin{proof}
  To help us prove this theorem we shall use three facts about compact right-topological semigroups:

  \begin{description}
    \item[(Fact 1)] 
      A subset $A$ of a semigroup is piecewise syndetic if and only if $\Closure{A} \cap K(\beta S) \ne \emptyset$.

    \item[(Fact 2)] 
      If $Y$ is a compact Hausdorff right-topological semigroup, $X \subseteq Y$ is a closed subsemigroup, and $K(Y) \cap X \ne \emptyset$, then $K(X) = K(Y) \cap X$.

    \item[(Fact 3)] 
          Given a family of discrete semigroups $\la S_i \ra_{i \in I}$, if we give $Y = \CartesianProduct_{ i \in I} \beta S_i$ the product topology and coordinatewise multiplication, then $Y$ is a compact Hausdorff right-topological semigroup and for every $\vec{x} \in \CartesianProduct_{i \in I} S_i$ we have that $p \mapsto \vec{x} \cdot p$ (for all $p \in Y$) is continuous, and $K(Y) = \CartesianProduct_{i \in I} K(\beta S_i)$.
  \end{description}

  \textbf{(Fact 3)} is a combination of \cite[\OutsideTheoremRef{Theorems 2.22 and 2.23 }]{Hindman:1998}, \textbf{(Fact 2)} follows from \cite[\OutsideTheoremRef{Theorems 1.65 and 2.7}]{Hindman:1998}, and \textbf{(Fact 1)} is \cite[\OutsideTheoremRef{Theorem 4.40}]{Hindman:1998}.
  With these preliminary results properly stated we now proceed to prove our theorem.

  Let $F \in \allNonemptyFiniteSubsetsOf{\calT}$, put $k = |F|$, and enumerate $F$ as $\{ f_1, f_2, \ldots, f_k \}$.
  Put $Y = \CartesianProduct_{t = 1}^k \beta S$.
  By \textbf{(Fact 3)}, giving $Y$ the product topology and pointwise multiplication, we have that $Y$ is a compact Hausdorff right-topological semigroup.
  
  We now define two special closed susbets of $Y$.
  To this end for every positive integer $N$ put
  \begin{align*}
    I_N = \Bigl\{ \bigl( x(m, a, t, f_1), &x(m, a, t, f_2), \ldots, x(m, a, t, f_k) \bigr) \Bigm\suchThat \\
      &\hspace{2em} \mbox{$m \in \bbN$, $a \in S^{m+1}$, and $t \in \UniformSplitting_m$ with $t_1 > N$} \Bigr\}
  \end{align*}
  and put $E_N = I_N \cup \{ (a, a, \ldots, a) \suchThat a \in S \}$.
  Let $E = \bigcap_{N = 1}^\infty c\ell_Y E_N$ and $I = \bigcap_{N = 1}^\infty c\ell_Y I_N$.
  Observe that both $E$ and $I$ are nonempty closed subsets of $Y$.
  We further claim one more important property on the relationship between $E$ and $I$:
  \begin{itemize}
    \item[(4)]
      $E$ is a subsemigroup of $Y$ and $I$ is an ideal of $E$.
  \end{itemize}

  For now, we will temporarily assume that property (4) is true and prove our theorem under the additional assumption of (4).
  (We will prove property (4) in the penultimate paragraph of this proof.)
  Since $A$ is piecwise syndetic, by \textbf{(Fact 1)} we can pick $p \in \Closure{A} \cap K(\beta S)$.
  Put $\overline{p} = (p, p, \ldots, p)$ and observe, by \textbf{(Fact 3)}, that $\overline{p} \in K(Y)$.
  Now, our immediate aim is to show that $\overline{p} \in E$.
  (So that $E \cap K(Y) \ne \emptyset$ and we can apply \textbf{(Fact 2)} to conclude $K(E) = K(Y) \cap E$.)
  Let $U$ be a neighborhood of $\overline{p}$ and pick $B_1, B_2, \ldots, B_k \in p$ such that $\CartesianProduct_{t = 1}^k \Closure{B_t} \subseteq U$.
  Pick $a \in \bigcap_{t=1}^k B_t$, then $(a, a, \ldots, a) \in U \cap E_N$ for every positive integer $N$.
  Hence $\overline{p} \in E$. and so $\overline{p} \in K(Y) \cap E$.
  Again, by \textbf{(Fact 2)} we conclude that $K(E) = K(Y) \cap E$.
  Since $K(E)$ is the smallest ideal in $E$ and by property (4) we have $\overline{p} \in K(E) \subseteq I \subseteq E$.
  Therefore $I_N \cap \CartesianProduct_{t = 1}^k A \ne \emptyset$ for ever positive integer $N$.
  Hence, by definition of $I_N$, we have that $A$ is a {\Jset}, under the additional assumption of statement (4).

  To finish our proof we now show property (4).
  Let $p$ and $q$ be elements of $E$, let $U$ be an open neighborhood of $p \cdot q$, and let $N$ be a positive integer.
  Since $r \mapsto  r \cdot q$ is continuous, by \textbf{(Fact 3)}, we can pick $V$ a neighborhood of $p$ such that $V \cdot q \subseteq U$.
  If $p \in I$, then we can pick $\vec{x} \in I_N \cap V$; otherwise pick $\vec{x} \in E_N \cap V$.
  If $\vec{x} \in I_N$, then pick a positive integer $m$, $a \in S^{m+1}$, and $t \in \UniformSplitting_m$ with $t_1 > N$ such that $\vec{x} = \bigl( x(m, a, t, f_1), x(m, a, t, f_2), \ldots, x(m, a, t, f_k) \bigr)$.
  In this case, put $M = t_m$; otherwise put $M = N$.

  Since $r \mapsto \vec{x} \cdot r$ is continuous, by \textbf{(Fact 3)}, we can pick $W$ a neighborhood of $q$ such that $\vec{x} W \subseteq U$.
  If $q \in I$, then pick $\vec{y} \in I_M \cap W$; otherwise pick $y \in E_M \cap W$.
  Then $\vec{x}\vec{y} \in E_N \cap U$, and if $p \in I$ or $q \in I$, then $\vec{x}\vec{y} \in I_N \cap U$.
  Hence it follows that $E$ is a subsemigroup of $Y$ and $I$ is an ideal of $E$.
  This completes the proof of this theorem.
\end{proof}
\begin{Remark}
  \label{remark:piecewise-syndetic-implies-Jset}
  The proof of \TheoremRef{piecewise-syndetic-implies-Jset} is modeled on the algebraic proof of van der Waerden's theorem found in \cite[\OutsideTheoremRef{Theorem 14.1}]{Hindman:1998}.
  The ideal for the algebraic proof of van der Waerden's theorem was based on Furstenberg and Katznelson ideas in \cite{Furstenberg:1989:CompactSemigrps}.
\end{Remark}

Therefore by \TheoremRef{piecewise-syndetic-implies-Jset} and \textbf{(Fact 1)} it follows that $K( \beta S) \subseteq J(S)$; which gives another proof of \CorollaryRef{J(S)-contains-smallest-ideal}.
We need \CorollaryRef{J(S)-contains-smallest-ideal} since we shall be performing a recursive construction on central sets (recall \DefinitionRef{central-set}) using the fact that a central set is a {\Jset}.

\section{Simpler Central Sets Theorem}
\label{section:simpler-central-sets-theorem}
With these preliminaries out of the way for {\Jsets} we now focus on deriving the new and simpler central sets theorem.
\begin{Definition}
  \label{definition:simpler-cset}
  Let $(S, \cdot)$ be a semigroup and $A \subseteq S$.
  We call $A$ a \Define{\Cset} if and only if there exist $m \colon \allNonemptyFiniteSubsetsOf{\calT} \to \bbN$, $\alpha \in \CartesianProduct_{F \in \allNonemptyFiniteSubsetsOf{\calT}} S^{m(F) + 1}$, and $\tau \in \CartesianProduct_{F \in \allNonemptyFiniteSubsetsOf{\calT}} \UniformSplitting_{m(F)}$ such that two conditions are satisfied:
  \begin{itemize}
    \item[(1)]
      if $F$ and $G$ are both elements of $\allNonemptyFiniteSubsetsOf{\calT}$ with $F \subsetneq G$, then $\tau(F)\bigl( m(F) \bigr) < \tau(G)(1)$, and

    \item[(2)]
      whenever $m$ is a positive integer, $G_1, G_2, \ldots, G_m \in \allNonemptyFiniteSubsetsOf{\calT}$ with $G_1 \subsetneq G_2 \subsetneq \cdots \subsetneq G_m$, and for each $i \in \{ 1, 2, \ldots, m \}$ we have $f_i \in G_i$, then 
      \[
        \prod_{i=1}^m x\bigl(m(G_i), \alpha(G_i), \tau(G_i), f_i \bigr) \in A.
      \]
  \end{itemize}
\end{Definition}

Our goal in this section is to prove the central sets theorem.
\begin{CentralSetsTheorem}
  \label{theorem:central-sets-theorem}
  Every central set in a semigroup is a {\Cset}.
\end{CentralSetsTheorem}

We now modify the sufficiency proof of \cite[\OutsideTheoremRef{Theorem 3.8}]{De:2008:NewStrongerCST} to show that every member of an idempotent in $J(S)$ is  a {\Cset}.
\begin{Theorem}
  \label{theorem:idempotents-and-Csets}
  Let $(S, \cdot)$ be a semigroup.
  If $p \in J(S)$ is an idempotent, that is, $p = p \cdot p$, then every $A \in p$ is a {\Cset}.
\end{Theorem}
\begin{proof}
  Since $A \in p$ and $p$ is an idempotent we have, by \cite[\OutsideTheoremRef{Lemma 4.14}]{Hindman:1998}, that $x^{-1}A^\star \in p$ for all $x \in A^\star$ where $A^\star = \{ x \in A \suchThat x^{-1}A \in p\}$. 
  We will recursively define our functions $m$, $\alpha$, and $\tau$ by the size of $F \in \FiniteSubsets(\calT)$ such that for $F \in \FiniteSubsets(\calT)$, the following hypotheses are satisfied:
  \begin{itemize}
    \item[(1)] If $\emptyset \ne G \subsetneq F$, then $\tau(G)\bigl( m(G) \bigr) < \tau(F)(1)$.
    
    \item[(2)] If $n \in \bbN$, $\emptyset \ne G_1 \subsetneq G_2 \subsetneq \cdots \subsetneq G_n = F$, and $\la f_i \ra_{i=1}^n \in \CartesianProduct_{i=1}^n G_i$, then \[ \prod_{i=1}^n x\bigl(m(G_i), \alpha(G_i), \tau(G_i), f_i)\bigr) \in A^\star.\]
  \end{itemize}

  Let $F \in \FiniteSubsets(\calT)$.
  First, assume that $|F| = 1$, that is $F = \{f\}$ for some sequence $f \in \calT$.
  Since $A^\star$ is a $J$-set, pick $m(F) \in \bbN$, $\alpha(F) \in S^{m(F)+1}$, and $\tau(F) \in \UniformSplitting_{m(F)}$ such that $x(m(F), \alpha(F), \tau(F), f) \in A^\star$. 

  Now assume that $|F| > 1$ and for all $\emptyset \ne G \subsetneq F$ we have defined $m(G)$, $\alpha(G)$, and $\tau(G)$ so that hypotheses (1) and (2) hold.
  Put 
  \begin{align*}
    M &= \bigl\{ \textstyle \prod_{i=1}^n x(m(G_i), \alpha(G_i), \tau(G_i), f_i) \;\bigm| \mbox{$n \in \bbN$, $\emptyset \ne G_1 \subsetneq G_2 \subsetneq \cdots \subsetneq G_n \subsetneq F$} \\
    &\hspace{12em} \mbox{and $\la f_i \ra_{i=1}^n \in \CartesianProduct_{i=1}^n G_i$} \bigr\}.
  \end{align*}
  Observe that since $F$ is finite, $M$ is also finite.
  By hypothesis (2) we have that $M \subseteq A^\star$.
  Put $B = A^\star \cap \bigcap_{x \in M} x^{-1}A^\star$, then $B \in p$ and so $B$ is a {\Jset}. 

  For each $\emptyset \ne G \subsetneq F$, put $l(G) = \tau(G)\bigl( m(G) \bigr)$ and put $k = \max\{ l(G) \;\suchThat \emptyset \ne G \subsetneq F \}$.
  By \LemmaRef{starting-position}, pick $m(F) \in \bbN$, $\alpha(F) \in S^{m(F)+1}$, and $\tau(F) \in \UniformSplitting_{m(F)}$ such that $\tau(F)(1) > k$ and for every $f \in F$, $x(m(F), \alpha(F), \tau(F), f) \in B$.
  Hypothesis (1) is satisfied, since $\tau(F)(1) > k \ge \max l(G)$ for all $\emptyset \ne G \subsetneq F$.
  We show that hypothesis (2) is also satisfied.
  Let $n \in \bbN$, $\emptyset \ne G_1 \subsetneq G_2 \subsetneq \cdots \subsetneq G_{n-1} \subsetneq G_n = F$, and $\la f_i \ra_{i=1}^n \CartesianProduct_{i=1}^n G_i$. 
  If $n = 1$, then $x(m(G_1), \alpha(G_1), \tau(G_1), f) \in B \subseteq A^\star$.
  Now assume that $n > 1$ and put $y = \prod_{i=1}^{n-1} x(m(G_i), \alpha(G_i), \tau(G_i), f_i)$. 
  By definition $y \in M$ and since $x(m(G_n), \alpha(G_n), \tau(G_n), f_n) = x(m(F), \alpha(F), \tau(F), f_n) \in B \subseteq y^{-1}A^\star$ we have
  \[
    \textstyle
    \prod_{i=1}^n x(m(G_i), \alpha(G_i), \tau(G_i), f_i) = y \cdot x(m(F), \alpha(F), \tau(F), f_n) \in A^\star.
  \]
  Hence hypotheses (1) and (2) are satisfied and this completes the proof for this theorem.
\end{proof}
\begin{Corollary}[Central Sets Theorem]
  \label{corollary:central-sets-theorem}
  In a semigroup every central set is a {\Cset}.
\end{Corollary}
\begin{proof}
  Let $(S, \cdot)$ be a semigroup and let $A \subseteq S$ be central.
  By \DefinitionRef{central-set} there exists an idempotent $p \in K(\beta S)$ with $A \in p$.
  By \CorollaryRef{J(S)-contains-smallest-ideal} $K(\beta S) \subseteq J(S)$ and so $A$ is also a {\Cset} by \TheoremRef{idempotents-and-Csets}.
\end{proof}

I'll finish this section by showing that, in fact {\Csets} can always be found as a member of an idempotent in $J(S)$.
To do this we shall need one more powerful algebraic lemma, which we won't prove in this paper.

\begin{Lemma}
  \label{lemma:super-lemma}     % Come up with a better name.
  Let $F$ be a set, $(D, \le)$ a directed set, and let $(S, \cdot)$ be a semigroup.
  Let $\la T_i \ra_{i \in D}$ be a decreasing family of subsets of $S$ such that for each $i \in D$ and $x \in T_i$, there exists $j \in D$ with $xT_j \subseteq T_i$.
  Put $\mathbf{Q} = \bigcap_{i \in D} c\ell_{\beta S}(T_i)$.
  Then $\mathbf{Q}$ is a compact subsemigroup of $\beta S$.
  Let $\la E_i \ra_{i \in D}$ and $\la I_i \ra_{i \in D}$ be decreasing families of nonempty subsets of $\CartesianProduct_{f \in F} S$ with the following properties:
  \begin{itemize}
    \item[(a)]
      For each $i \in D$, $I_i \subseteq E_i \subseteq \CartesianProduct_{f \in F} S$. 

    \item[(b)]
      For each $i \in D$ and every $\vec{x} \in I_i$, there exists $j \in D$ such that $\vec{x}E_j \subseteq I_i$.

    \item[(c)]
      For each $i \in D$ and every $\vec{x} \in E_i \setminus I_i$, there exists $j \in D$ such that $\vec{x} E_j \subseteq E_i$ and $\vec{x}I_j \subseteq I_i$.
  \end{itemize}

  Let $Y = \CartesianProduct_{f \in F} \beta S$, let $E = \bigcap_{i \in D} c\ell_{Y}(E_i)$, and let $I = \bigcap_{i \in D} c\ell_Y(I_i)$.
  Then $E$ is a subsemigroup of $\CartesianProduct_{f \in F} \mathbf{Q}$ and $I$ is an ideal of $E$. 
  Additionally, if either 
  \begin{itemize}
    \item[(d)]
      for each $i \in D$, $T_i = S$ and $\{ a \in S \;\suchThat \overline{a} \not\in E_i \}$ is \emph{not} piecewise syndetic, or

    \item[(e)]
      for each $i \in D$ and each $a \in T_i$, $\overline{a} \in E_i$,
  \end{itemize}
then given any $p \in K(\mathbf{Q})$, we have $\overline{p} \in E \cap K(\CartesianProduct_{f \in F} \mathbf{Q}) = K(E) \subseteq I$. 
\end{Lemma}
\begin{proof}
  This is proved in \cite[\OutsideTheoremRef{Lemma 14.9}]{Hindman:1998}.
\end{proof}

We now modify the necessity proof of \cite[\OutsideTheoremRef{Theorem 3.8}]{De:2008:NewStrongerCST} to show that a {\Cset} can always be found as a member of an idempotent in $J(S)$.
\begin{Theorem}
  \label{theorem:idempotents-in-J(S)}
  Let $(S, \cdot)$ be a semigroup and let $A \subseteq S$.
  If $A$ is a {\Cset}, then there exists an idempotent $p \in J(S)$ with $A \in p$.
\end{Theorem}
\begin{proof}
  Pick $m \colon \allNonemptyFiniteSubsetsOf{\calT} \to \bbN$, $\alpha \in \CartesianProduct_{F \in \allNonemptyFiniteSubsetsOf{\calT}} S^{m(F)+1}$, and $\tau \in \CartesianProduct_{F \in \allNonemptyFiniteSubsetsOf{\calT}} \UniformSplitting_{m(F)}$ as guaranteed by the definition of a {\Cset}. 
  For each $F \in \allNonemptyFiniteSubsetsOf{\calT}$ define 
  \begin{align*}
    T(F) &= \bigl\{ \textstyle \prod_{i=1}^n x(m(F_i), \alpha(F_i), \tau(F_i), f_i) \;\bigm\suchThat \mbox{$n \in \bbN$, for each $i \in \{1, 2, \ldots, n\}$, }\\
 &\hspace{8em}\mbox{$F_i \in \allNonemptyFiniteSubsetsOf{\calT}$, $F \subsetneq F_1 \subsetneq F_2 \subsetneq \cdots \subsetneq F_n$, $\la f_i \ra_{i=1}^n \in \CartesianProduct_{i=1}^n F_i$} \bigr\}.
  \end{align*}

  Observe that for each $F \in \allNonemptyFiniteSubsetsOf{\calT}$, we have $T(F)$ is a nonempty subset of $A$, and the collection $\{ T(F) \;\suchThat F \in \allNonemptyFiniteSubsetsOf{\calT} \}$ has the finite intersection property since $T(F \cup G) \subseteq T(F) \cap T(G)$ for all $F$, $G \in \allNonemptyFiniteSubsetsOf{\calT}$.
  Therefore $\mathbf{Q} = \bigcap_{F \in \allNonemptyFiniteSubsetsOf{\calT}} c\ell_{\beta S}\bigl(T(F)\bigr)$ is a closed nonempty subset of $\beta S$. 

  We show that $\mathbf{Q}$ is in fact a subsemigroup of $\beta S$.
  To see that $\mathbf{Q}$ is a subsemigroup it suffices, by \cite[\OutsideTheoremRef{Theorem 4.20}]{Hindman:1998}, to show that for all $F \in \allNonemptyFiniteSubsetsOf{\calT}$ and for every $y \in T(F)$, there exists $G \in \allNonemptyFiniteSubsetsOf{\calT}$ such that $yT(G) \subseteq T(F)$. 
  So let $F \in \allNonemptyFiniteSubsetsOf{\calT}$ and $y \in T(F)$.
  Pick $n \in \bbN$, for every $i \in \{1, 2, \ldots, n\}$ pick $F_i \in \allNonemptyFiniteSubsetsOf{\calT}$ with $F \subsetneq F_1 \subsetneq F_2 \subsetneq \cdots \subsetneq F_n$, and $\la f_i \ra_{i=1}^n \in \CartesianProduct_{i=1}^n F_i$ such that $y = \prod_{i=1}^n x(m(F_i), \alpha(F_i), \tau(F_i), f_i)$. 
  We show that $yT(F_n) \subseteq T(F)$.
  Let $z \in T(F_n)$ and pick $m \in \bbN$, for every $i \in \{1, 2, \ldots, m\}$ pick $G_i \in \allNonemptyFiniteSubsetsOf{\calT}$ with $F_n \subsetneq G_1 \subsetneq G_2 \subsetneq \cdots \subsetneq G_m$, and pick $\la g_i \ra_{i=1}^m \in \CartesianProduct_{i=1}^m G_i$ such that $z = \prod_{i=1}^m x(m(G_i), \alpha(G_i), \tau(G_i), g_i)$.

  For each $i \in \{1, 2, \ldots, n+m\}$ define 
  \[
    H_i = 
    \begin{cases}
      F_i & \mbox{if $i \in \{1, 2, \ldots, n\}$,} \\
      G_{i-n} & \mbox{if $i \in \{n+1, n+2, \ldots, n+m\}$,}
    \end{cases}
  \]
  and define the sequence
  \[
    h_i =
    \begin{cases}
      f_i & \mbox{if $i \in \{1, 2, \ldots, n\}$,} \\
      g_{i-n} & \mbox{if $i \in \{n+1, n+2, \ldots, n+m\}$.}
    \end{cases}
  \]
  Then $F \subsetneq H_1 \subsetneq H_2 \subsetneq \cdots \subsetneq H_{n+m}$, $\la h_i \ra_{i=1}^{n+m} \in \CartesianProduct_{i=1}^{n+m} H_i$, and 
  \[
    \textstyle
    yz = \prod_{i=1}^{n+m} x(m(H_i), \alpha(H_i), \tau(H_i), h_i) \in T(F)
  \]
  Hence $\mathbf{Q}$ is a subsemigroup of $\beta S$. 

  We now claim that $K(\mathbf{Q}) \subseteq \overline{A} \cap J(S)$. 
  If this claim is true, then we are done since any idempotent in $K(\mathbf{Q})$ will establish the theorem. 
  % Since $\mathbf{Q} \subseteq c\ell_{\beta S} (A)$, it evident that $K(\mathbf{Q}) \subseteq c\ell_{\beta S} (A)$. 

  Now let $p \in K(\mathbf{Q})$ and $B \in p$.
  Using \LemmaRef{super-lemma} we shall show that $B$ is a $J$-set. 
  Let $F \in \FiniteSubsets(\calT)$ and put $D = \{G \in \FiniteSubsets(\calT) \; \suchThat F \subseteq G \}$. 
  Observe $\mathbf{Q} = \bigcap_{G \in D} c\ell_{\beta S} (T(G))$.

  For $G \in D$, define $I(G) \subseteq \CartesianProduct_{f \in F} S$ as follows: for $w \in \CartesianProduct_{f \in F} S$, $w \in I(G)$ if and only if there is some $n \in \bbN \setminus \{1\}$ such that 
  \begin{itemize}
    \item[(1)]
      there exist disjoint nonempty sets $C_1$ and $C_2$ with $\{1, 2, \ldots, n\} = C_1 \cup C_2$, 
      
    \item[(2)]
      there exists a strictly increasing sequence $\la G_i \ra_{i=1}^n$ in $\FiniteSubsets(\calT)$ with $G \subsetneq G_1$, and

    \item[(3)]
      there exists $\sigma \in \CartesianProduct_{i \in C_1} G_i$, such that for every $f \in F$, if $\gamma_f \in \CartesianProduct_{i=1}^n G_i$ is defined by 
  \[
    \gamma_f(i) = 
    \begin{cases}
      \sigma(i) & \mbox{if $i \in C_1$,} \\
      f & \mbox{if $i \in C_2$,} 
    \end{cases}
  \]
then $w(f) = \prod_{i=1}^n x(m(G_i), \alpha(G_i), \tau(G_i), \gamma_f(i))$.
  \end{itemize}
  For each $G \in D$, define $E(G) = I(G) \cup \{ \overline{b} \;\suchThat b \in T(G) \}$. 

  We claim that $\la E(G) \ra_{G \in D}$ and $\la I(G) \ra_{G \in D}$ satisfy statements \textsl{(a)}, \textsl{(b)}, \textsl{(c)}, and \textsl{(e)} of \LemmaRef{super-lemma}.

  Assume, temporarily, that our claim is true. 
  Then by \LemmaRef{super-lemma} if $Y = \CartesianProduct_{f \in F} \beta S$, $E = \bigcap_{G \in D} c\ell_Y(E(G))$, and $I = \bigcap_{G \in D} c\ell_Y(I(G))$, then $E$ is a subsemigroup of $Y$, $I$ is an ideal of $E$, and for every $p \in K(\mathbf{Q})$, $\overline{p} = (p, p, \ldots, p) \in K(E) \subseteq I$. 
  Since $\CartesianProduct_{f \in F} c\ell_{\beta S} (B)$ is a neighborhood of $\overline{p}$ we can pick $w \in I(F) \cap \CartesianProduct_{f \in F} c\ell_{\beta S} (B)$. 
  Pick $n \in \bbN\setminus \{1\}$, $C_1$, $C_2$, $\la G_i \ra_{i=1}^n$, and $\sigma \in \CartesianProduct_{i \in C_1} G_i$ as guaranteed by the definition of $I(F)$. 
  Put $r = |C_2|$ and enumerate $C_2$ as a strictly increasing sequence $h_1$, $h_2$, \dots, $h_r$.
  Put $u = \sum_{i=1}^r m(G_i)$.
  We define $c \in S^{u+1}$ and $t \in \UniformSplitting_u$ such that for all $f \in F$, $x(u, c, t, f) \in B$.

  Define $c \in S^{u+1}$ as follows:
  \[
    c(1) = 
    \begin{cases}
      \alpha(G_1)(1) & \mbox{if $h_1 = 1$,} \\
      \bigl(\prod_{i=1}^{h_1-1} x(m(G_i), \alpha(G_i), \tau(G_i), \sigma(i))\bigr) \cdot \alpha(G_{h_1})(1) & \mbox{if $h_1 > 1$.}
    \end{cases}
  \]
  For each positive integer $j$ with $1 < j < m(G_{h_1})$ put $c(j) = \alpha(G_{h_2})(j)$, and for each positive integer $j$ with $1 \le j \le m(G_{h_1})$ define $t_j = \tau(G_{h_1})(1)$. 

  Now for each $s \in \{1, 2, \ldots, u-1\}$ put $v_s = \sum_{i=1}^s m(G_{h_i})$, and define
  \[
    c(v_s+1) = 
    \begin{cases}
      \alpha(G_{h_s})(m(G_{h_s}+1))\alpha(G_{h_{s+1}})(1) & \mbox{if $h_{s+1} = h_s + 1$,} \\

      \alpha(G_{h_s})(m(G_{h_s}+1))\cdot & \\
      \hspace{2em}\bigl(\prod_{i=h_{s}+1}^{h_{s+1}-1} x(m(G_i), \alpha(G_i), \tau(G_i), \sigma(i))\bigr)\alpha(G_{h_{s+1}})(1) & \mbox{if $h_{s+1} > h_s+1$.}
    \end{cases}
  \]
  For each $s \in \{1, 2, \ldots, u-1\}$ and every positive integer $j$ with $v_s < j \le \sum_{i=1}^{s+1} m(G_{h_i})$ put $t_j = \tau(G_{h_{s+1}})(j-u)$. 
  Finally, we define
  \[
    c(u+1) = 
    \begin{cases}
      \alpha(G_{h_r})(m(G_n) + 1) & \mbox{if $h_r = n$,} \\
      \alpha(G_{h_r})(m(G_{h_r}+1))\prod_{i=h_r+1}^n x(m(G_i), \alpha(G_i), \tau(G_i), \sigma(i)) & \mbox{if $h_r < n$.}
    \end{cases}
  \]
  Then for every $f \in F$, $x(u, c, t, f) \in B$, and so $B$ is a $J$-set.

  We now prove our claim that the families $\la E(G) \ra_{G \in D}$ and $\la I(G) \ra_{G \in D}$ satisfy statements \textsl{(a)}, \textsl{(b)}, \textsl{(c)}, and \textsl{(e)} of \LemmaRef{super-lemma}.

  By definition of $\la E(G) \ra_{G \in D}$ and $\la I(G) \ra_{G \in D}$ it is immediate that statements \textsl{(a)} and \textsl{(e)} are satisfied. 

  We now show statement (b) which states that for every $G \in D$ and each $w \in I(G)$ there exists $H \in D$ such that $w E(H) \subseteq I(G)$.
  Let $G \in D$ and $w \in I(G)$. 
  Pick $n \in \bbN$, $C_1$, $C_2$, $\la G_i \ra_{i=1}^n$, and $\sigma \in \CartesianProduct_{i \in C_1} G_i$ as guaranteed by the definition of $I(G)$..
  We show that $w E(G_n) \subseteq I(G)$. 
  Let $z \in E(G_n)$. 
  First assume that $z = \overline{b}$ for some $b \in T(G_n)$. 
  Pick $m \in \bbN$, for each $i \in \{1, 2, \ldots, m\}$ pick $F_i \in \FiniteSubsets(\calT)$ with $G_n \subsetneq F_1 \subsetneq F_2 \subsetneq \cdots \subsetneq F_n$, and pick $\la f_i \ra_{i=1}^n \in \CartesianProduct_{i=1}^n F_i$ such that $b = \prod_{i=1}^m x(m(F_i), \alpha(F_i), \tau(F_i), f_i)$. 
  Put $D_1 = C_1 \cup \{n+1, n+2, \ldots, n+m\}$ and for each $i \in \{1, 2, \ldots, n+m\}$ put 
  \[
    H_i = 
    \begin{cases}
      G_i & \mbox{if $i \le n$,} \\
      F_{i-n} & \mbox{if $i > n$.}
    \end{cases}
  \]
  Define $\rho \in \CartesianProduct_{i \in D_1} H_i$ by 
  \[
    \rho(i) =
    \begin{cases}
      \sigma(i) & \mbox{if $i \le n$,} \\
      f_{i-n} & \mbox{if $i > n$.}
    \end{cases}
  \]
  Then with $n+m$, $D_1$, $C_2$, $\la H_i \ra_{i=1}^{n+m}$, and $\rho$ we have that $w \cdot z \in I(G)$. 

  Now we assume that $z \in I(G_n)$. 
  Pick $m \in \bbN$, $D_1$, $D_2$, $\la F_i \ra_{i=1}^m$, and $\rho$ as guaranteed by the definition of $I(G_n)$. 
  Put $E_1 = C_1 \cup \{n + i \;\suchThat i \in D_1 \}$ and put $E_2 = C_2 \cup \{ n + i \;\suchThat i \in D_2 \}$.
  For each $i \in \{1, 2, \ldots, n+m\}$ put
  \[
    H_i = 
    \begin{cases}
      G_i & \mbox{if $i \le n$,} \\
      F_{i-n} & \mbox{if $i > n$.}
    \end{cases}
  \]
  Define $\mu \in \CartesianProduct_{i \in E_i} H_i$ by
  \[
    \mu(i) = 
    \begin{cases}
      \sigma(i) & \mbox{if $i \le n$,} \\
      \rho(i) & \mbox{if $i > n$.}
    \end{cases}
  \]
  Then with $n+m$, $E_1$, $E_2$, $\la H_i \ra_{i=1}^{n+m}$, and $\mu$ we have that $w \cdot z \in I(G)$.

  We now verify statement \textsl{(c)} which states that for all $G \in D$ and every $w \in E(G) \setminus I(G)$, there exists $H \in D$ such that $wE(H) \subseteq E(G)$ and $wI(H) \subseteq I(G)$.
  So let $G \in D$ and $w \in E(G) \setminus I(G)$. 
  Pick $b \in T(G)$ such that $w = \overline{b}$. 
  Pick $n \in \bbN$, for each $i \in \{1, 2, \ldots, n\}$ pick $G_i \in \FiniteSubsets(\calT)$ with $G \subsetneq G_1 \subsetneq G_2 \subsetneq \cdots \subsetneq G_n$, and pick $\la f_i \ra_{i=1}^n \in \CartesianProduct_{i=1}^n G_i$ such that $b = \prod_{i=1}^n x(m(G_i), \alpha(G_i), \tau(G_i), f_i)$. 
  Then similar to what we have done above, we have that $wE(G_n) \subseteq E(G)$ and $wI(G_n) \subseteq I(G)$.
\end{proof}

\section{Equivalence between two versions of the Central Sets Theorem}
\label{section:equivalence-of-central-set-theorems}
In this last section we show that our ``new'' definition of a {\Jset} is equivalent to the ``old'' definition of a {\Jset}.
(I point out again that the main advantage of the definition of a {\Jset} given here is that it is easier to work with.)
To show that the simpler version of the noncommutative central sets theorem \CorollaryRef{central-sets-theorem} is equivalent to the older formulation \cite[\OutsideTheoremRef{Corollary 3.10}]{De:2008:NewStrongerCST} it is necessary and sufficient to show that the simpler formulation of a {\Jset}, \DefinitionRef{new-Jset}, is equivalent to the older formulation of a {\Jset} in \cite[\OutsideDefinitionRef{Definition 3.3(e)}]{De:2008:NewStrongerCST}.

To help us prove this equivalence we use the following ``rewriting'' lemma.
\begin{Lemma}
  \label{lemma:rewriting-lemma}
  Let $(S, \cdot)$ be a semigroup, let $n$ be a positive integer, let $c \in S^{n+1}$, let $H \in \NonUniformSplitting_n$, and $F \in \allNonemptyFiniteSubsetsOf{\calT}$.
  Fix $b \in S$ and for each $f \in F$ define $g_f \in \calT$ by $g_f(t) = f(t) \cdot b$ for every $t \in \bbN$.
  Then there exists $m \in \bbN$, $a \in S^{n+1}$, and $t \in \UniformSplitting_m$ such that for all $f \in F$ we have
  \[
    x(m, a, t, f) = \Bigl(\prod_{i = 1}^n \bigl(c(i) \cdot \prod_{u \in H(i)} g_f(u) \bigr) \Bigr) c(m + 1).
  \]
\end{Lemma}
\begin{proof}
  Put $H(0) = \emptyset$ and for each $s \in \{ 0, 1, \ldots, n \}$ define $h_s = \sum_{i = 0}^s |H(i)|$.
  Put $m = h_n$ and enumerate $\bigcup_{i = 1}^n H(i)$ as a strictly increasing sequence $t_1 < t_2 < \cdots < t_m$.

  We shall adopt some temporary terminology: for every $f \in F$ we say \[\textstyle\Bigl(\prod_{i = 1}^n \bigl( c(i) \cdot \prod_{u \in H(i)} g_f(u) \bigr)\Bigr) \cdot c(n+1)\] has \Define{proper representation} if and only if $\Bigl(\prod_{i = 1}^n \bigl( c(i) \cdot \prod_{u \in H(i)} g_f(u) \bigr)\Bigr) \cdot c(n+1) = x(m, a, t, f)$ where $a \in S^{m + 1}$ is defined by
  \[
    a(j) =
    \begin{cases}
      c(1) & \mbox{if $j = 1$,}\\
      b    & \mbox{if $s \in \{ 0, 1, \ldots, m - 1 \}$ and $2 + h_s \le j \le h_{s + 1}$, and} \\
      b \cdot c(s + 1) & \mbox{if $s \in \{1, 2, \ldots, m\}$ and $j = 1 + h_s.$}
    \end{cases}
  \]
  
  We prove, by induction on $n$, that $\Bigl(\prod_{i = 1}^n \bigl( c(i) \cdot \prod_{u \in H(i)} g_f(u) \bigr)\Bigr) \cdot c(n+1)$ has proper representation for every $f \in F$.

  First, suppose $n = 1$, then $\Bigl(\prod_{i = 1}^1 \bigl( c(i) \cdot \prod_{u \in H(i)} g_f(u) \bigr)\Bigr) \cdot c(2) = c(1) \cdot f(t_1) \cdot  b \cdot f(t_2) \cdot b \cdots f(t_m) \cdot b \cdot c(2)$.
  In this case $h_0 = 0$, $h_1 = m$, and so $s$ can only be 0 or 1.
  If $2 + h_0 = 2 \le j \le h_1 = n$, then by definition of $a$ we have $a(j) = b$ for all $j \in \{ 2, 3, \ldots, m \}$.
  Also since $1 + h_1 = n + 1$, we have $c( n + 1) = b \cdot a(2)$.
  Therefore
  \begin{align*}
    c(1) \cdot f(t_1) \cdot  b \cdot f(t_2) \cdot b \cdots f(t_m) \cdot b \cdot c(2) &= a(1) \cdot f(t_1) \cdot a(2) \cdot f(t_2) \cdot a(3) \cdots \\
    &\hspace{2em}a(m) \cdot f(t_m) \cdot a(m + 1).
  \end{align*}
  Hence $\Bigl(\prod_{i = 1}^1 \bigl( c(i) \cdot \prod_{u \in H(i)} g_f(u) \bigr)\Bigr) \cdot c(2)$ has proper representation.

  Now let $n > 1$ and assume that $\Bigl(\prod_{i = 1}^{n - 1} \bigl( c(i) \cdot \prod_{u \in H(i)} g_f(u) \bigr)\Bigr) \cdot c(n)$ has proper representation, say with $\Bigl(\prod_{i = 1}^{n - 1} \bigl( c(i) \cdot \prod_{u \in H(i)} g_f(u) \bigr)\Bigr) \cdot c(n) = x(m, a, t, f)$.
  Then we have
  \begin{align*}
    \Bigl(\prod_{i = 1}^{n} \bigl( c(i) \cdot \prod_{u \in H(i)} g_f(u) \bigr)\Bigr) \cdot c(n + 1) &= \Bigl(\prod_{i = 1}^{n - 1} \bigl( c(i) \cdot \prod_{u \in H(i)} g_f(u) \bigr)\Bigr) \cdot c(n) \\
    &\hspace{5em}\cdot \bigl(\prod_{t \in H(n)} g_f(t)\bigr) \cdot c(n + 1), \\
    &= x(m, a, t, f) \cdot \bigl(\prod_{t \in H(n)} g_f(t)\bigr) \cdot c(n + 1).
  \end{align*}
  Now by our base case $a(m + 1) \bigl( \prod_{t \in H(n)} g_f(t) \bigr) c(n + 1)$ has proper representation, say with
  \[
    a(m + 1) \bigl( \prod_{t \in H(n)} g_f(t) \bigr) c(n + 1) = x(p, d, u, f).
  \]
  By properly translating the indices for $u$ and $d$ it follows that \[\textstyle\Bigl(\prod_{i = 1}^{n} \bigl( c(i) \cdot \prod_{u \in H(i)} g_f(u) \bigr)\Bigr) \cdot c(n + 1)\] has proper representation.
\end{proof}

\begin{Theorem}
  \label{theorem:equivalence-between-notions-of-J-sets}
  Let $(S, \cdot)$ be a semigroup and $A \subseteq S$.
  The following two statements are equivalent:
  \begin{itemize}
    \item[(a)] $A$ is a {\Jset}.

    \item[(b)] For all $F \in \allNonemptyFiniteSubsetsOf{\calT}$, there exists $m \in \bbN$, $a \in S^{m+1}$, and $H \in \calI_m$ such that for all $f \in F$ we have $\prod_{i=1}^m\bigl(a(i) \prod_{t \in H(i)} f(t) \bigr)a(m+1) \in A$.
  \end{itemize}
\end{Theorem}
\begin{proof}
  Assume that $A$ is a {\Jset}.
  Then $A$ obviously satisfies statement (b) by putting $H_i = \{ t_i \}$ for each $i \in \{ 1, 2, \ldots , m\}$.

  Now assume $A$ satisfies the statement (b).
  Let $F \in \allNonemptyFiniteSubsetsOf{\calT}$ and fix $b \in S$.
  For each $f \in F$ define $g_f \in \calT$ by $g_f(t) = f(t) \cdot b$.
  By assumption we can pick $n \in \bbN$, $c \in S^{n + 1}$, and $H \in \NonUniformSplitting_n$ such that for every $f \in F$ we have $\Bigl(\prod_{i = 1}^n\bigl( c(i) \cdot \prod_{t \in H(i)} g_f(t)\bigr)\Bigr) \cdot a(m + 1) \in A$.
  By \LemmaRef{rewriting-lemma} there exist $m \in \bbN$, $a \in S^{m + 1}$, and $t \in \UniformSplitting_m$ such that for every $f \in F$ we have $x(m, a, t, f) = \Bigl(\prod_{i = 1}^n\bigl( c(i) \cdot \prod_{t \in H(i)} g_f(t)\bigr)\Bigr) \cdot a(m + 1)$.
\end{proof}

\emergencystretch=1em
\printbibliography
\end{document}